\documentclass{article}

\begin{document}

\title{The Ding functional, Berndtsson convexity and moment maps}
%

\maketitle
\newcommand{\bC}{{\bf C}}
\newcommand{\bR}{{\bf R}}
\newcommand{\bP}{{\bf P}}
\newtheorem{lem}{Lemma}
\newtheorem{thm}{Theorem}
\newtheorem{cor}{Corollary}
\newtheorem{prop}{Proposition}

\

{\it \ \ \ \ \ \ \ \ \ \ \ \ \ \ \  To Jean-Michel Bismut, for his 65th. birthday}

\

\section{Introduction}
This article is largely expository in nature. The main results we discuss are not new but our goal is to fit them into a framework which does not seem to have been considered explicitly before and which we hope has some conceptual value. 

As background, we  recall briefly a standard package of ideas which arises when a Lie group $G$ acts by holomorphic isometries on a K\"ahler manifold $Z$ and we have an equivariant moment map
\begin{equation} \mu:Z\rightarrow {\rm Lie}(G)^{*}. \end{equation}
Initially, we have in mind the case when $G$ is compact. Then there is a complexified group $G^{c}$ which acts holomorphically on $Z$. One central topic is then to relate the {\it symplectic quotient} $\mu^{-1}(0)/G$ and the complex quotient $Z^{s}/G^{c}$, where $Z^{s}\subset Z$ is an appropriate subset of \lq\lq stable points''. This amounts to solving an equation
$\mu(g(z_{0}))=0$ for $g\in G^{c}$ when $z_{0}$ is a stable point. This problem can be cast in variational shape. Given $z_{0}$, we define a $1$-form $\theta$ on the group $G^{c}$  by 
\begin{equation}  \theta_{g}(\delta g)= 2 \langle \mu(g(z_{0})), i  \delta g\ g^{-1}\rangle \end{equation}
(where $\delta g$ is a tangent vector to $G^{c}$ at $g$). One finds that this is a closed $1$-form and can be expressed as the derivative of a function ${\cal F}$, unique up to  constant. Moreover, everything is invariant under the group $G$, so ${\cal F}$ can be regarded as a function on the symmetric space $G^{c}/G$. An important point is that the function is {\it convex} along geodesics in this space. Solving the moment map equation is equivalent to minimising ${\cal F}$. (In Section 4 below we will take the geometry of this general situation a bit further, in the discussion of evolution equations.)

Beginning with Atiyah and Bott \cite{kn:AB}, it has been realised that many problems in differential geometry can be fitted into this framework, with an infinite dimensional group $G$. In other words we can cast a differential geometric equation of interest in the shape $\mu=0$ where $\mu$ is the moment map for a  group action. In this article we want to discuss the case of K\"ahler metrics, specifically we consider the following two differential geometric problems.
\begin{enumerate} \item In a fixed K\"ahler class on any compact K\"ahler manifold $X$, study metrics of {\it constant scalar curvature}.
\item Study {\it K\"ahler-Einstein metrics} on a Fano manifold (i.e. when $c_{1}(X)$ is a K\"ahler class). 
\end{enumerate}
Beginning with Fujiki \cite{kn:Fujiki} (and earlier work of Quillen and others in the case of one complex dimension) it has been understood that the first problem fits neatly into this conceptual form. We assume here some familiarity with this picture, as in \cite{kn:Don1}, \cite{kn:Don2}. The space $G^{c}/G$ in question is interpreted as the space of K\"ahler metrics in a given class and the functional ${\cal F}$ is the Mabuchi functional. The convexity of this functional has been exploited to very good effect in the study of constant scalar curvature metrics by Chen \cite{kn:Chen} and many others.

From one point of view, the second problem is a special case of the first. Clearly a K\"ahler-Einstein metric has constant scalar curvature and conversely a simple integration argument shows that in the class $c_{1}(X)$ a constant scalar curvature metric is K\"ahler-Einstein. However this point of view is not completely satisfactory since the K\"ahler-Einstein problem is fundamentally simpler than the constant scalar curvature one (for example, the equation  can be reduced to a second order PDE, rather than fourth order) and one would like to have a set-up tailored specifically to it.  Further, there is a natural functional in the literature on this problem, the {\it Ding functional} \cite{kn:Ding} and this was shown to be convex in important work of Berndtsson \cite{kn:BB}. The problem we address in this article is to fit the Ding functional and its convexity into the standard conceptual package outlined above.  To achieve this we will use the same group as in the constant scalar curvature theory but use a different metric on the infinite dimensional space on which it acts, exploiting the special features of the situation. In this point of view Berndtsson's calculation amounts to the statement that the metric is positive definite. 

The author is very grateful to Dietmar Salamon, Song Sun and the referee for helpful comments.

\section{Preliminaries}

We review standard material concerning the homogeneous space $Sp(n,\bR)/U(n)$.
Let $U$ be a real vector space of dimension $2n$ and let $c_{n}=1$ if $n$ is even and $i$ if $n$ is odd. For $\alpha,\beta\in \Lambda^{n}(U)\otimes \bC$ write
\begin{equation} \langle \alpha,\beta \rangle = c_{n} \alpha\wedge \overline{\beta}. \end{equation}
This is a $\Lambda^{2n}\otimes \bC$-valued indefinite Hermitian form on $\Lambda^{n}(U)\otimes \bC$.
Now suppose that $U^{*}$ has a symplectic structure $\omega\in \Lambda^{2}(U)$. The set $M$ of compatible complex structures on $U$ can be identified with an open subset of the
Lagrangian Grassmanian of $U\otimes \bC$ (i.e. $n$-dimensional complex subspaces which are Lagrangian with respect to the complex extension of $\omega$). Thus we can identify $M$ with
a subset of $\bP(\Lambda^{n} U\otimes \bC)$. Let $N\subset \Lambda^{n} U\otimes \bC$ be the corresponding cone, with origin deleted. In the standard model, with $U=\bC^{n}$, the group $Sp(n,\bR)$ acts transitively on $M$ and the stabiliser of the standard complex structure is $U(n)$. The points of $N$ are just those which are $\bR^{+}\times Sp(n,R)$ equivalent to
\begin{equation}\alpha_{0}= dz_{1}\dots dz_{n}. \end{equation}

With the normalisation we adopted in (3),  and using the symplectic volume form $\omega^{n}/n!$  to identify $\Lambda^{2n}(U)$ with $\bR$, we have $\langle\alpha,\alpha\rangle>0$ for $\alpha\in N $.
 The tangent space $TN$ of $N$ at $\alpha$ is a complex linear subspace of $\Lambda^{n}(U)\otimes \bC$, containing $\alpha$. The basic observation is that the Hermitian form is {\it negative definite} on the orthogonal complement of $\alpha$.  Equivalently, writing $\beta$ for a tangent vector, the expression
  $$   \frac{\langle \beta, \beta\rangle}{\langle \alpha, \alpha\rangle}
  -\left\vert\frac{\langle \alpha,\beta\rangle}{\langle \alpha,\alpha\rangle}\right\vert^{2}$$
  
   defines a negative-definite metric on the tangent space of $M$ at $[\alpha]$. Up to a normalising factor, this is the negative of the standard homogeneous, K\"ahler, structure on $M$. All this is clear from the fact that $Sp(n,\bR)$ acts transitively on $M$ and the stabiliser $U(n)$ of the standard point acts irreducibly on the tangent space. To see things more explicitly work at the standard model point $\alpha_{0}$.
The tangent space to the decomposable forms at $\alpha_{0}$ is the direct sum $\Lambda^{n,0}\oplus \Lambda^{n-1,1}$. Using the isomorphism
$$ \Lambda^{n-1,1}\otimes (\Lambda^{n,0})^{*} = \Lambda^{0,1} \otimes U^{*}$$
and the symplectic form, we get an isomorphism
   $$  \Lambda^{n-1,1}\rightarrow  U^{*}\otimes U^{*}\otimes \Lambda^{n,0}$$
   and the orthogonal complement of $\alpha_{0}$ in the tangent space of $N$
corresponds to $s^{2}(U^{*})\otimes \Lambda^{n,0}$. However we do not need to use this below.

  From another point of view, if we regard $N$ as a $\bC^{*}$-bundle over $M$ we are saying that the Hermitian structure on this bundle defined by $\langle \alpha, \alpha \rangle$ generates the K\"ahler form as its curvature. From yet another point of view, the K\"ahler structure on $M$ is the K\"ahler quotient by the action of $S^{1}$ on the indefinite linear structure restricted to $N$. 

In the next section we will need the following identity. 
Let $v\in U^{*}\otimes \bC$ and $\alpha \in N$. The contraction $v\dashv \alpha$ is in
$\Lambda^{n-1}(U)$. We have then 
\begin{equation}  c_{n} \ \omega \wedge (v\dashv \alpha) \wedge (\overline{i v\dashv\alpha})= (-1)^{n-1} \vert v'\vert^{2} \langle \alpha,\alpha \rangle, \end{equation}
where $v=v'+v''$ is the decomposition of $v$ defined by the  complex structure corresponding to $\alpha$ and $\vert \ \vert$ is the standard metric defined by $\alpha$ and $\omega$. To see this we can suppose that $\alpha$ is the standard structure in (4), $\omega=\frac{i}{2}\sum_{a} dz_{a} d\overline{z}_{a}$ and write 
$$   v= \sum v'_{a} \frac{\partial}{\partial z_{a}} + \sum v''_{a} \frac{\partial}{\partial\overline{z}_{a}}, $$
from which point the formula (5) is an easy calculation.

\section{K\"ahler metrics on Fano manifolds}
Let $(X,\omega)$ be a compact symplectic $2n$-manifold with $-i\omega$ the curvature form of a connection $\nabla$ on a Hermitian line bundle $L$.
We consider the natural bundle $\underline{M}\rightarrow X$ with fibre $M$.  Thus sections of $\underline{M}$ are almost-complex structures on $X$ compatible with the symplectic form. We write ${\cal J}$ for the space of sections. The K\"ahler structure on $M$ induces a formal K\"ahler structure on ${\cal J}$ in a straightforward way. (In general, if we have a fibre bundle $E\rightarrow B$ where $B$ has a volume form and the fibres have K\"ahler structures then there is a formal K\"ahler structure on the space of sections.) Let ${\cal G}$ be the group of automorphisms of  $L$, as a line bundle with connection,
covering symplectomorphisms of $X$. This acts on ${\cal J}$ preserving the structures and the moment map can be identified as the \lq\lq hermitian scalar curvature''. This is the theory developed in \cite{kn:Don1}. Most often one restricts attention to the space ${\cal J}_{int}\subset {\cal J}$ of integrable almost-complex structures, so one is describing K\"ahler metrics and in that case the Hermitian scalar curvature becomes the ordinary scalar curvature. To fit in with the discussion in the preceding section we should normalise the moment map by subtracting the average value of the scalar curvature.

Now we begin our  different approach. Consider the $L$-valued forms $\Omega^{n}(L)$. For such forms $\alpha, \beta$ we define
$\langle \alpha, \beta\rangle$ as before except that we also use the Hermitian form on $L$. Then $\langle \alpha, \beta\rangle$ is a $2n$-form on $X$. Write $d_{\nabla}$ for the coupled exterior derivative on $L$-valued forms. So we have an integration-by-parts formula
\begin{equation}   \int_{X} \langle d_{\nabla}\sigma, \beta\rangle= (-1)^{n}\int_{X} \langle \sigma , d_{\nabla}\beta \rangle, \end{equation}
for $\sigma \in \Omega^{n-1}(L), \beta\in \Omega^{n}(L)$,  where we extend the definition of $\langle\ ,\ \rangle$ in the obvious way. By definition, the group ${\cal G}$ acts on $L$ and on $X$ and so on $\Omega^{n}(L)$. Recall the basic fact from symplectic geometry that the Lie algebra of ${\cal G}$ is $C^{\infty}(X)$. That is, given a Hamiltonian function $H$ we get a Hamiltonian vector field $v_{H}$ on $X$ but also a lift to $L$ in which we rotate the fibres by $iH$.  This leads to a formula for the infinitesimal action of the Lie algebra $C^{\infty}(X)$ on $\Omega^{n}(L)$ as
\begin{equation}  i R_{\alpha}(H) = v_{H}(\dashv  d_{\nabla}\alpha) - d_{\nabla}( v_{H}\dashv \alpha) - i H \alpha. \end{equation} 
(The factor $i$ is included for convenience later.)

Restrict to the case when the first Chern class of $L$ is equal to the first Chern class of any compatible almost-complex structure. Let $\underline{N}$ be the bundle over $X$ with fibre $N$, coupled to $L$ as above. So sections of $\underline{N}$ are certain $L$-valued $n$- forms on $X$ and the Chern class condition above is just the condition that global sections of $\underline{N}$ exist. The space $\underline{N}$ is a principal $\bC^{*}$-bundle over $\underline{M}$. 
\begin{lem}
Any section $\alpha$ of $\underline{N}$ with $d_{\nabla} \alpha=0$ projects to an integrable almost complex structure. Conversely if $J\in {\cal J}_{int}$ is an integrable complex structure it has a lift to a section $\alpha$ of $\underline{N}$ with $d_{\nabla} \alpha=0$ and $\alpha$ is unique up to multiplication by a constant in $\bC^{*}$. 
\end{lem}
In one direction, suppose that $\alpha$ is a section of $\underline{N}$ with $d_{\nabla}\alpha=0$. Choose a local trivialisation of $L$ in the neighbourhood of a point of $X$. Thus in this trivialisation $\alpha$ can be regarded as a complex $n$-form on $X$ and
\begin{equation}d\alpha = \theta \alpha\end{equation} for some $1$-form $\theta$. This is one version of the standard integrability condition (The almost-complex structure defines a decomposition of the complex tangent space $TX\otimes \bC= T'X\oplus T''X$ and (8) implies that $\Gamma(T''X)$ is closed under Lie bracket.) In the other direction, suppose that $J$ is an integrable almost-complex structure. Then $L$ becomes a positive holomorphic line bundle over the complex manifold $X$ and our hypothesis implies that $c_{1}(X)>0$. By a standard vanishing theorem we have $H^{0,1}(X)=0$ and so $K^{-1}_{X}$ is isomorphic to $L$ as a holomorphic line bundle. Fixing such an isomorphism, we have a canonical element $\alpha$ given in local complex co-ordinates by $s \otimes dz_{1} \dots dz_{n}$ where 
$$  s= \frac{\partial}{\partial z_{1}} \wedge \dots \wedge \frac{\partial}{\partial z_{n}}.$$
Then the fact that $\nabla s$ has type $(1,0)$ (since $s$ is holomorphic) implies that
$  d_{\nabla}\alpha=0$. If also $d_{\nabla}(f \alpha)=0$ then $df\wedge \alpha=0$ which implies that $\overline{\partial} f=0$ and so $f$ is a constant (since $X$ is compact).

 Let $\hat{{\cal J}}_{int}$ be the sections $\alpha$ of $\underline{\hat{M}}$ with $d_{\nabla} \alpha=0$.
Thus the lemma assets that $\hat{{\cal J}}_{int}$ is a $\bC^{*}$-bundle over ${\cal J}_{int}$.
The space $\hat{{\cal J}}_{int}$ is a subspace of the vector space $\Omega^{n}(L)$ which has a Hermitian form given by
\begin{equation}  \langle \langle \alpha, \beta \rangle \rangle = \int_{X} \langle \alpha, \beta \rangle. \end{equation}
Clearly $  \langle \langle \alpha, \alpha \rangle \rangle >0$ for any $\alpha\in \hat{{\cal J}}_{int}$.

 Let $T_{\alpha}\subset \Omega^{n}(L)$consist of the forms $\beta$ with $d_{\nabla}\beta=0$ and which at each point $x\in X$ lie in the corresponding vertical tangent space of $\underline{N}$. This can be viewed as the tangent space of
$\hat{{\cal J}}_{int} $ at $\alpha$. (There is a technical point that $\hat{{\cal J}}_{int}, {\cal J}_{int}$ could be singular spaces so we are working with the Zariski tangent space, but we do not need to go into this here.) Clearly  $\alpha$ lies in $T_{\alpha}$. The main result we need is an analogue of the standard discussion in Section 2.
\begin{thm}
The Hermitian form $\langle\langle\ ,\ \rangle\rangle$ is negative-definite on the orthogonal complement of $\alpha$ in $T_{\alpha}$.
\end{thm}
As in Section 2, what this means is that any $\beta \in T_{\alpha}$ we have
\begin{equation}   \langle \langle\beta, \beta\rangle\rangle - \frac{\vert \langle\langle \alpha, \beta \rangle\rangle\vert^{2} }{ \langle\langle \alpha, \alpha \rangle\rangle} \leq 0\end{equation}
with equality if and only if $\beta$ is a  complex multiple of $\alpha$. 
To prove the Theorem we first extend the operator $R_{\alpha}$ by linearity to complex valued functions on $X$. Thus if $f=f_{1}+ i f_{2}$ with $f_{1}, f_{2}$ real we define
 $$ R_{\alpha} f  =  i d_{\nabla} (v_{f}\dashv \alpha) \ -  f \alpha $$
where $v_{f}= v_{f_{1}} + i v_{f_{2}}$.  Then $T_{\alpha}$ contains the image of $R_{\alpha}$ since it is a complex vector space and, by naturality of our constructions, it contains $R_{\alpha} H$ for $H$ real. We first establish the inequality (10) for $\beta$ of the form $R_{\alpha} f$. To do this, define $Pf$ by the pointwise orthogonal decomposition
\begin{equation}    R_{\alpha} f = (P f) \alpha + \alpha^{\perp}, \end{equation}
where $\langle\alpha^{\perp},\alpha\rangle$ vanishes at each point of $X$. In other words
\begin{equation}   (Pf)\langle \alpha, \alpha\rangle= \langle R_{\alpha} f, \alpha\rangle. \end{equation}
 Now we calculate $\langle\langle R_{\alpha} f, R_{\alpha} f\rangle\rangle$ in three different ways.
\begin{enumerate}
\item The integration by parts formula (6) implies that
$$  \langle\langle d_{\nabla}(v_{f}\dashv \alpha) , R_{\alpha}f\rangle\rangle =  0$$
since $d_{\nabla} R_{\alpha} f =0$. Thus
\begin{equation} \langle\langle R_{\alpha} f, R_{\alpha} f\rangle\rangle= -\langle\langle R_{\alpha} f,  f\alpha \rangle\rangle. \end{equation}
Now using the definition of $Pf$ we get
\begin{equation}   \langle\langle R_{\alpha} f, R_{\alpha} f\rangle\rangle= - \int_{X} \overline{f} P(f) \langle \alpha,\alpha \rangle \end{equation}
The same argument shows that
\begin{equation}   \langle \langle R_{\alpha} f , R_{\alpha} g \rangle \rangle= -\int_{X} \overline{g} P(f) \langle \alpha,\alpha \rangle, \end{equation}
which shows that $P$ is a self-adjoint operator with respect to the $L^{2}$ norm defined by the measure $\langle \alpha, \alpha \rangle$.

\item 
Applying the integration by parts formula again we get
\begin{equation}\langle\langle f\alpha, R_{\alpha} f\rangle\rangle= -\langle \langle f \alpha, f \alpha\rangle\rangle+ (-1)^{n} \langle\langle d_{\nabla} (f \alpha),  i v_{f} \dashv f \alpha \rangle \rangle \end{equation}
Using the definition of curvature, $d_{\nabla}^{2} = -i \omega\wedge$, the identity $d_{\nabla} R_{\alpha} f =0$ gives
\begin{equation}  - i\omega \wedge ( v_{f}\dashv \alpha)) - d_{\nabla}( i f \alpha) =0\end{equation}
(Of course this can be derived directly.)  So we get
\begin{equation}   \langle\langle R_{\alpha} f, R_{\alpha} f\rangle\rangle= \int_{X} \vert f \vert^{2}  \langle \alpha, \alpha \rangle +  c_{n}  \omega \wedge (v_{f} \dashv \alpha) \wedge (\overline{i v_{f}\dashv \alpha}). \end{equation}
Applying (5), we have
$$  c_{n} \omega \wedge (v_{f} \dashv \alpha) \wedge (\overline{i v_{f}\dashv \alpha})= \vert v'_{f}\vert^{2} \langle \alpha, \alpha \rangle. $$
Under the metric defined by $\alpha$ the component $v'_{f}$ corresponds to $\partial f $ so
$\vert v'_{f}\vert^{2} = \vert \partial f \vert^{2}$ and we get
\begin{equation} \langle \langle R_{\alpha} f,R_{\alpha} f\rangle\rangle =\int_{X} (\vert f \vert^{2} - \vert \partial f \vert^{2}) \langle \alpha, \alpha \rangle \end{equation}
\item
Going back to the formula $ R_{\alpha} f = P f \alpha +\alpha^{\perp}$, we have
$$  \langle R_{\alpha} f, R_{\alpha} f\rangle= \vert P f\vert^{2} \langle \alpha, \alpha\rangle
+ \langle \alpha^{\perp}, \alpha^{\perp}\rangle . $$
By the discussion in Section 1 $\langle \alpha^{\perp},\alpha^{\perp}\rangle \leq 0$ with equality if and only if $\alpha^{\perp}=0$.  Thus
\begin{equation} \langle\langle R_{\alpha} f, R_{\alpha} f \rangle \rangle\leq \int_{X} \vert Pf\vert^{2} \langle\alpha, \alpha \rangle, \end{equation}
with strict inequality unless $R_{\alpha} f= (Pf)\alpha$.

\end{enumerate}

Suppose $\lambda$ is an eigenvalue of $P$. Then (14) and (19) imply that $\lambda\leq 1$ and if $\lambda=1$ the eigenspace consists of the constant functions. On the other hand (14) and (20) imply that $\lambda^{2} \geq \lambda$, so if $\lambda\neq 1$ we have $\lambda\leq 0$. Moreover if the case $\lambda=0$ occurs the corresponding eigenspace is the kernel of $R_{\alpha}$. 
Write $(\ ,\ )$ for the $L^{2}$ inner product on complex valued functions on $X$ corresponding to the measure $\langle \alpha,\alpha \rangle$. We have then
$$   (Pf,f) \leq \frac{\vert(f,1)\vert^{2}}{(1,1)}, $$
with equality if and only if $Pf=0$. Using (15) this shows that (10) holds for $\beta$ in the image of $R_{\alpha}$. 

To complete the proof of Theorem 1 it suffices to show that (10) holds if $\beta\in T_{\alpha} $ is orthogonal, under $\langle\langle \ ,\ \rangle \rangle$, to the image of $R_{\alpha}$. But by the integration by parts formula
$$   \langle \langle \beta, R_{\alpha} f \rangle \rangle = -\int_{X} \overline{f} \langle \beta, \alpha \rangle, $$
and if this  vanishes for all $f$ we must have $\langle \beta,\alpha\rangle=0$ at each point of $X$. It follows immediately from Section 2 that $\langle \langle \beta,\beta \rangle \rangle \leq 0$ with equality if and only if $\beta=0$.

Theorem 1 means that the indefinite form -$\langle\langle\ , \ \rangle\rangle$
induces a positive definite K\"ahler metric on ${\cal J}_{int}$. (The fact
that it is K\"ahler can be seen  by regarding it as a symplectic quotient,
for example.) This is different from the metric considered at the beginning of this section and as we will see below, it achieves our goal of placing the K\"ahler-Einstein theory in the general package outlined at the beginning of the paper. But we should emphasise that Theorem 1 is essentially a restatement of Berndtsson's convexity theorem (in fact a special case of that theorem) from a different point of view, and the proof we have given is probably essentially the same as the standard one but in a different notation.

\section{Moment maps and K\"ahler-Ricci flow}

 The group ${\cal G}$ acts linearly on the symplectic vector space $\Omega^{n}(L)$ and the moment map is simply given by:
$$ \alpha \mapsto -\frac{1}{2} \left( H \mapsto {\rm Re} \langle\langle \alpha, R_{\alpha} H\rangle\rangle\right). $$
When restricted to forms $\alpha$ with $d_{\nabla} \alpha=0$ we have $$\langle \langle \alpha, R_{\alpha} H\rangle\rangle= \langle\langle \alpha, H \alpha \rangle\rangle $$
(using the integration-by-parts formula on one term). Thus the moment map, $\mu_{0}$, for the action
on ${\hat{\cal J}}_{int}$ is just
$$  \alpha \mapsto -\frac{1}{2}\left(H \mapsto \int_{X} H \langle \alpha, \alpha\rangle\right). $$
Or in other words we have $\mu_{0}(\alpha)=-\frac{1}{2}\langle \alpha, \alpha \rangle$, using the pairing between $2n$-forms and functions. Since the action of ${\cal G}$ commutes with the $S^{1}$ action, the moment map on the symplectic quotient ${\cal J}_{int}$ is defined by
$$  \mu_{0}([\alpha])= -\frac{1}{2}\frac{\langle \alpha,\alpha\rangle}{\langle\langle \alpha,\alpha\rangle \rangle}. $$

As before we can adjust the moment map by a constant and it is more convenient to use
$$\mu([\alpha])= -\frac{1}{2} \left(\frac{\langle \alpha,\alpha\rangle}{\langle\langle
\alpha,\alpha\rangle \rangle} - C \omega^{n}\right) , $$
where $C^{-1}$ is the integral of $\omega^{n}$ over $X$. With this normalisation $\mu([\alpha])$ is a $2n$-form of integral $0$. 

For any $[\alpha]\in {\cal J}_{int}$ the pair $(\omega, [\alpha])$ defines a K\"ahler metric $g(\omega, [\alpha])$. Then we have
\begin{prop}
For $[\alpha]\in {\cal J}_{int}$ we have $\mu([\alpha])=0$ if and only if $g(\omega,[\alpha])$ is K\"ahler-Einstein.
\end{prop}
This is just a matter of tracing through definitions. Recall first the standard formulation of the K\"ahler-Einstein condition on a Fano manifold, which we denote by $X_{J}$, with canonical bundle $K$. Let $h$ be a hermitian metric on $K^{-1}$. This can be viewed, algebraically, as a volume form $\Omega_{h}$.
In terms of local co-ordinates, if 
$$ \vert \frac{\partial}{\partial z_{1}}\wedge \dots\wedge \frac{\partial}{\partial z_{n}}\vert_{h}^{2} = V$$
then
$$   \Omega_{h} = V \left(\frac{i}{2}\right)^{n} dz_{1}d\overline{z}_{1}\dots dz_{n}d\overline{z}_{n} . $$
On the other hand, the metric $h$ defines a Chern connection on $K^{-1}$  with curvature form $-i \omega_{h}$ say, and we have another volume form
$\omega_{h}^{n}$. K\"ahler-Einstein metrics correspond to metrics $h$ such that $\omega_{h}>0$ and $\Omega_{h}=\omega_{h}^{n}$.

In our formulation, we  take $X_{J}$ to be the complex manifold corresponding to $[\alpha]$ and we use $\alpha$ to identify the line bundle $L$ with  $K^{-1}$. Since $L$ has a metric, this gives a hermitian metric $h$ on $K^{-1}$ and the connection $\nabla$ can now be viewed as the Chern connection on $K^{-1}$, so $\omega=\omega_{h}$. We also have $\Omega_{h}=\langle \alpha,\alpha\rangle$. The condition $\mu([\alpha])=0$ becomes 
$$    \Omega_{h}= C \langle\langle \alpha,\alpha\rangle\rangle \omega_{h}^{n}$$
and we merely need to match up normalisations by choosing  $\alpha$ so that
$$\langle\langle \alpha,\alpha\rangle\rangle= C^{-1}. $$

Next we turn to the functional associated to the problem. Recall that in the formal picture developed in \cite{kn:Mabuchi}, \cite{kn:Don1},\cite{kn:Don2} and elsewhere, one interprets an orbit of the complexification of ${\cal G}$ in ${\cal J}_{int}$ as the set of complex structures in a fixed isomorphism class. In other words, while the group ${\cal G}^{c}$ does not exist its orbits still make sense. Fix attention on one orbit corresponding to a fixed complex manifold structure  $X_{J}$ as above.  One interprets the space ${\cal G}^{c}/{\cal G}$ as the space ${\cal H}$ of metrics $h$ on $K^{-1}$ with $\omega_{h}>0$. If we fix a reference metric $h_{0}$ then ${\cal H}$ becomes   the space of K\"ahler potentials via $h=e^{-\phi} h_{0}$ and $\omega_{h}= \omega_{0}+ i \partial \overline{\partial} \phi$. Unravelling the definitions, one finds that the functional ${\cal F}$ has derivative:
$$  \delta {\cal F}=  \int (\delta \phi) (-C \omega_{h}^{n}+\frac{\Omega_{h}}{\int \Omega_{h}}). $$
This means that $$ {\cal F}= -I(h)- \log \int \Omega_{h}, $$
where the functional $I$ is defined by
$$   \delta I= C \int (\delta \phi)\omega_{h}^{n}. $$
This is the standard description of the Ding functional ${\cal F}$.

We want to finish by discussing the K\"ahler-Ricci flow in this picture. Return to the general conceptual situation considered at the beginning of the article  when we have a Lie group $G$ acting by holomorphic isometries on a K\"ahler manifold $(Z,\Omega)$ and an equivariant moment map
$\mu:Z\rightarrow {\rm Lie}(G)^{*}$. Suppose first that we have a invariant metric on ${\rm Lie}(G)$, then the function $E(z)= \frac{1}{2}\Vert \mu(z)\Vert^{2}$ is a smooth function on $Z$. The gradient flow of this function is
$$   \frac{dz}{dt}= I \rho_{z}(\mu(z))  , $$ where
$\rho_{z}: {\rm Lie}(G)\rightarrow TZ_{z}$ is the derivative of the action and we have used the metric to identify the Lie algebra with its dual. One interesting feature of this flow is that it preserves $G^{c}$ orbits and on each such orbit co-incides with the gradient flow of the functional ${\cal F}$, with respect to the homogenous metric on $G^{c}/G$. In fact one has
     $\frac{d \cal F}{dt} = - 2 E$ so one sees that $\frac{d^{2}{\cal F}}{dt^{2}}\geq 0$ 
The corresponding flows appear in a number of infinite dimensional situations, in particular in the case of K\"ahler metrics the Calabi flow fits into this framework with fixed points the extremal metrics.

 More generally, suppose given a function $B: {\rm Lie}(G)^{*}\rightarrow \bR$, invariant under the co-adjoint action. Then
$H=B(\mu)$ is a $G$-invariant function on $Z$ and its gradient flow is 
\begin{equation}  \frac{dz}{dt}= I \rho_{z} (DB(\mu(z)), \end{equation}
where $DB: {\rm Lie}(G)^{*}\rightarrow {\rm Lie}(G) $ is the derivative of $B$. Thus it is still true that the flow preserves $G^{c}$ orbits. Now we have
$$  \frac{d{\cal F}}{dt} = - (DB(\mu), \mu).$$
Suppose that $B$ is a non-negative convex function, with $B(0)=0$. Then $(DB(\mu), \mu)\geq B(\mu)$ so we have
\begin{equation} \frac{d{\cal F}}{dt} \leq - H. \end{equation}

In our situation, on a Fano manifold, with the action of the group ${\cal G}$ on ${\cal J}_{int}$ and the \lq\lq Berndtsson metric'' on the latter, we first identify functions with $2n$-forms using the form $\omega^{n}$. (This is just for convenience). So now we regard $\mu([\alpha])$ as the function
$$   \mu=    \frac{\langle \alpha,\alpha\rangle}{\omega^{n}\langle\langle \alpha,\alpha\rangle\rangle} - C. $$
(Notice that if we put $\theta=\mu+C$ then $\log \theta$ is the the \lq\lq Ricci potential'': $i\partial\overline{\partial} \log \theta = \rho-\omega$, where $\rho$ is the Ricci form, with the normalisation  that the integral of $\theta \omega^{n}$ is $1$.)
Define 
$$  B(\mu)= \int_{X} \left[ (\mu+C)\log (\frac{\mu}{C}+1) -\mu\right] \omega^{n}. $$
(This does not strictly fit into the picture above since it is not defined on the whole of ${\rm Lie}({\cal G})^{*}$, but it is defined on the image of $\mu$ which is clearly all that is needed.) It is clear that this is a convex nonnegative function, vanishing only at $\mu=0$.  We have $DB(\mu) = \log (\frac{\mu}{C}+1)$ and the flow (18) translates  to 
$$   \frac{d\omega}{dt} =  -i\partial\overline{\partial} \log \frac{\Omega_{h}}{\omega_{h}^{n}}. $$
which is the K\"ahler-Ricci flow. This functional $H=B(\mu)$ is the functional introduced by  He in \cite{kn:He}. In terms of a metric $h$ on $K^{-1}$, normalised so that the integral of $\Omega_{h}$ is equal to the integral of $\omega_{h}^{n}$, we have 
$$   H= C \int_{X} \log\left(\frac{\Omega_{h}}{\omega_{h}^{n}}\right)\Omega_{h}. $$

The general statements above  give that, under the K\"ahler-Ricci flow:
\begin{enumerate}
\item $H$ is decreasing;
\item the Ding functional ${\cal F}$ satisfies $\frac{d{\cal F}}{dt} \leq -H$.
\end{enumerate}

The first property was found by  He (equation (2.4) in \cite{kn:He}). The second has perhaps not been noticed before.
Of course these are simple calculations which can be checked without  the conceptual picture above. On the other hand the conceptual picture is helpful in suggesting what functionals to consider. There is a useful corollary of these inequalities.

\begin{cor}
Suppose $\omega_{t}$ is a K\"ahler-Ricci flow for $0\leq t<\infty$. Suppose there is a sequence of times $ t(\nu)\rightarrow \infty$ such that $\omega_{t(\nu)}$ converges to a metric which is not K\"ahler-Einstein. Then the Ding functional ${\cal F}(\omega_{t})$ tends to $-\infty$ as $t\rightarrow \infty$.
\end{cor}
 This follows immediately from the statements above and the fact that the zeros of $H$ are exactly the K\"ahler-Einstein metrics. 

The \lq\lq moment map'' picture may also be useful in the study of K\"ahler-Ricci solitons, along the lines of Szekelyhidi's work on extremal metrics \cite{kn:Sz}.

\end{document}